\documentclass{amsart}

\DeclareMathSymbol{\twoheadrightarrow}  {\mathrel}{AMSa}{"10}

\def\Z{{\mathbb Z}}

\def\RR{{\mathbb R}}
\def\F{{\mathbb F}}
\def\P{{\mathbb P}}

\def\Sn{{\mathbf S}_n}
\def\An{{\mathbf A}_n}
\def\RR{{\mathfrak R}}
\def\Perm{\mathrm{Perm}}
\def\Gal{\mathrm{Gal}}
\def\Pic{\mathrm{Pic}}

\def\End{\mathrm{End}}
\def\Aut{\mathrm{Aut}}
\def\div{\mathrm{div}}

\def\fchar{\mathrm{char}}

\def\A{\mathbf{A}}

\newtheorem{thm}{Theorem}[section]
\newtheorem{lem}[thm]{Lemma}
\newtheorem{cor}[thm]{Corollary}
\newtheorem{prop}[thm]{Proposition}
\theoremstyle{definition}

\title[Del Pezzo surfaces and jacobians]
{Del Pezzo surfaces of degree $2$ and  jacobians without complex
multiplication}

\author[Yu.\ G.\ Zarhin]{Yu.\ G.\ Zarhin}
\address{Department of Mathematics, Pennsylvania State University,
University Park, PA 16802, USA}

\address{Institute for Mathematical Problems in Biology, Russian Academy of
Sciences, Pushchino, Moscow Region, Russia}
\email{zarhin\char`\@math.psu.edu} \dedicatory{To my friend Sergei
Vostokov}

\begin{document}

\maketitle

\section{Notations and Statements}
In a series of his articles
\cite{ZarhinMRL,ZarhinMRL2,ZarhinTexel,ZarhinMMJ,ZarhinBSMF} the
author constructed explicitly $m$-dimensional abelian varieties
without non-trivial endomorphisms for every $m>1$. This
construction may be described as follows. Let $K_a$ be an
algebraic closure of a perfect field $K$ with $\fchar(K)\ne 2$.
Let $n=2m+1$ or $2m+2$. Let us choose an $n$-element set $\RR \in
K_a$ that constitutes a Galois orbit over $K$ and assume, in
addition, that the Galois group of $K(\RR)$ over $K$ is ``big"
say, coincides with full symmetric group $\Sn$ or the alternating
group $\An$. Let $f(x) \in K[x]$ be the irreducible polynomial of
degree $n$, whose set of roots coincides with $\RR$. Let us
consider the hyperelliptic curve $C_f:y^2=f(x)$ over $K_a$ and let
$J(C_f)$ be its jacobian which is the $m$-dimensional abelian
variety. Then the ring $\End(J(C_f))$ of all $K_a$-endomorphisms
of $J(C_f)$ coincides with $\Z$ if either $n>6$ or $\fchar(K)\ne
3$.

The aim of this paper is to construct abelian threefolds without
complex multiplication, using jacobians of non-hyperelliptic
curves of genus $3$. It is well-known that these curves are smooth
plane quartics and closely related to Del Pezzo surfaces of degree
$2$. (We refer to \cite{Manin,
Iskovskikh,IskovskikhSh,Demazure,Dolgachev, Dolgachev2,T} for
geometric and arithmetic  properties of Del Dezzo surfaces. In
particular,   relations between the degree $2$ case and plane
quartics are discussed in detail in  \cite{Demazure,Dolgachev,
Dolgachev2}). On the other hand, Del Pezzo surfaces of degree $2$
could be obtained by blowing up seven points on the projective
plane $\P^2$ when these points are in {\sl general position},
i.e., no three points lie on a one line, no six on a one conic
(\cite[\S 3]{Iskovskikh}, \cite[Th. 1 on p. 27]{Demazure}).

In order to describe our construction,  let us start with the
projective plane $\P^2$ with homogeneous coordinates $(x:y:z)$.
Let us consider a $7$-element set $B\subset \P^2(K_a)$ of points
in general position and assume that the absolute Galois group
$\Gal(K)$ of $K$ permutes elements of $B$ in such a way that $B$
constitutes a Galois orbit. We write $Q_B$ for the $6$-dimensional
$\F_2$-vector space of maps $\varphi:B \to \F_2$ with $\sum_{b\in
B}\varphi(b)=0$. The action of $\Gal(K)$ on $B$ provides $Q_B$
with the natural structure of $\Gal(K)$-module. Let $G_B$ be the
image of $\Gal(K)$ in the group $\Perm(B)\cong \mathbf{S}_7$ of all
permutations of $B$. Clearly, $Q_B$ carries a natural structure of
faithful $\Perm(B)$-module and the structure homomorphism $\Gal(K)
\to \Aut(Q_B)$ coincides with the composition of
$\Gal(K)\twoheadrightarrow G_B$ and
$G_B\subset\Perm(B)\hookrightarrow \Aut(Q_B)$.

Let $H_B$ be the $K_a$-vector space of homogeneous cubic forms in
 $x,y,z$ that vanish on $B$. It follows
from  proposition 4.3  and corollary 4.4(i) in Ch. 5, \S 4 of
\cite{Har} that $H_B$ is $3$-dimensional and $B$ coincides with
the set of common zeros of elements of $H_B$. Since $B$ is
$\Gal(K)$-invariant, $H_B$ is defined over $K$, i.e., it has a
$K_a$-basis $u,v,w$ such that the forms $u,v,w$ have coefficients
in $K$.

We write $V(B)$ for the Del Pezzo surface of degree $2$ obtained
by blowing up $B$. Then $V(B)$ is a smooth projective surface that
is defined over $K$ (see Remark 19.5 on pp. 89--90 of
\cite{Manin}). We write
$$g_B:V(B)\to \P^2$$
for the corresponding birational map defined over $K$. Recall that
for each $b\in B$ its preimage $E_b$ is a a smooth projective
rational curve with self-intersection number $-1$. By definition,
$g_B$ establishes a $K$-biregular isomorphism between
$V(B)\setminus \bigcup_{b\in B}E_b$ and $\P^2\setminus B$.
Clearly,
$$\sigma(E_b)=E_{\sigma(b)}\quad \forall \ b\in B,
\sigma\in\Gal(K).$$ Let $\Omega_{V(B)}$ be the canonical
(invertible) sheaf on $V(B)$. Let us consider the line $L:z=0$ as
a divisor in $\P^2$. Clearly, $B$ does not meet the $K$-line $L$;
otherwise, the whole $\Gal(K)$-orbit $B$ lies in $L$ which is not
true, since no $3$ points of $B$ lie on a one line. It is known
\cite[Sect. 25.1 and 25.1.2 on pp. 126--127] {Manin} that
$$K_{V(B)}:=-3g_B^*(L)+\sum_{b\in B}E_b=-g_B^*(3L)+\sum_{b\in B}E_b$$ is a
{\sl canonical divisor} on $V(B)$. Clearly, for each form $q\in
H_B$ the rational function $\frac{q}{z^3}$ on $\P^2$ satisfies
$\div(\frac{q}{z^3})+3L\ge 0$, i.e., $\frac{q}{z^3}\in
\Gamma(\P^2,3L)$. Also $\frac{q}{z^3}$ is defined and vanishes at
every point of $B$. It follows easily that $\frac{q}{z^3}$ (viewed
as rational function on $V(B)$) lies in
$\Gamma(V(B),3g_B^*(L)-\sum_{b\in B}E_b)=\Gamma(V(B), -K_{V(B)})$.
Since $\Gamma(V(B), -K_{V(B)})$ is $3$-dimensional \cite[theorem
24.5 on p. 121]{Manin},
$$\Gamma(V(B), -K_{V(B)})=K_a \cdot \frac{u}{z^3}\oplus K_a \cdot \frac{v}{z^3}\oplus
K_a \cdot \frac{w}{z^3}\ .$$ Using proposition 4.3 in  \cite[Ch.
5, \S 4]{Har}, one may easily get a well-known fact that the
sections of $\Gamma(V(B), -K_{V(B)})$ have no common zeros on
$V(B)$. This gives us a regular anticanonical map
$$\pi:V(B)\stackrel{g_B}{\longrightarrow}\P^2
\stackrel{(u:v:w)}{\longrightarrow}\P^2$$ which is obviously
defined over $K$. It is  known  that $\pi$ is a regular double
cover map, whose ramification curve is a smooth quartic
$$C_B \subset \P^2$$
(see \cite[pp. 67--68]{Demazure}, \cite[Ch. 9]{Dolgachev}).
Clearly, $C_B$ is a genus $3$ curve defined over $K$. Let $J(B)$
be the jacobian of $C_B$; clearly, it is a three-dimensional
abelian variety defined over $K$. We write $\End(J(B))$ for the
ring of $K_a$-endomorphisms of $J(B)$.

The following assertion is based on  Lemmas 1-2 on pp. 161--162 of
\cite{Dolgachev}.

\begin{lem}
\label{points} Let $J(B)_2$ be the kernel of multiplication by $2$
in $J(B)(K_a)$. Then the Galois modules $J(B)_2$ and $Q_B$ are
canonically isomorphic.
\end{lem}

Using Lemma \ref{points} and results of
\cite{ZarhinMRL,ZarhinBSMF}, one may obtain the following
statement.

\begin{thm}
\label{endo} Let $B\subset \P^2(K_a)$ be a $7$-element set of
points in general position. Assume that $\Gal(K)$ permutes
elements of $B$ and
 the image of $\Gal(K)$ in $\Perm(B)\cong \mathbf{S}_7$ coincides either
with the full symmetric group $\mathbf{S}_7$ or with the alternating
group $\A_7$.  Then $\End(J(B))=\Z$.
\end{thm}

 This leads to a question: how to
construct such $B$ in general position? The next lemma provides us
with desired construction.

\begin{lem}
\label{general} Let $f(t) \in K[t]$ be a separable irreducible
degree $7$ polynomial, whose Galois group $\Gal(f)$ is either
$\mathbf{S}_7$ or $\A_7$. Let $\RR_f\subset K_a$ be the $7$-element set
of roots of $f$. Then the $7$-element set
$$B_f =\{(\alpha^3:\alpha:1)\mid \alpha\in\RR_f\}\subset
\P^2(K_a)$$ is in general position.
\end{lem}

Clearly, $\Gal(K)$ permutes transitively elements of $B_f$ and
 the image of $\Gal(K)$ in $\Perm(B)$ coincides either
with  $\mathbf{S}_7$ or with   $\A_7$; in particular, $B_f$ constitutes a
Galois orbit. This implies the following statement.

\begin{cor}
Let $f(t) \in K[t]$ be a separable irreducible degree $7$
polynomial, whose Galois group $\Gal(f)$ is either $\mathbf{S}_7$ or
$\A_7$. Then  $\End(J(B_f))=\Z$.
\end{cor}

\section{Proofs}
\begin{proof}[Proof of Lemma \ref{points}]
Let $\Pic(V(B))$ be the Picard group of $V(B)$ over $K_a$. It is
known \cite[Sect. 25.1 and 25.1.2 on pp. 126--127] {Manin} that
$\Pic(V(B))$ is a free commutative group of rank $8$ provided with
the natural structure of Galois module. More precisely, it has
canonical generators $l_0=$ the class of $g_B^*(L)$  and
$\{l_b\}_{b\in B}$ where $l_b$ is the class of the exceptional
curve $E_b$. Clearly, $l_0$ is Galois invariant and
$$\sigma(l_b)=l_{\sigma(b)} \quad \forall \ b\in B, \sigma\in
\Gal(K).$$ Clearly, the class of $K_{V(B)}$ equals
$-3l_0+\sum_{b\in B}l_b$ and obviously is Galois-invariant. There
is a non-degenerate Galois invariant symmetric intersection form
$$(, ): \Pic(V(B)) \times \Pic(V(B)) \to \Z.$$
In addition (ibid),
$$(l_0, l_0)=1, (l_b,l_0)=0, (l_b,l_b)=-1, (l_b,l_{b'})=0 \quad
\forall\ b \ne b'.$$ Clearly, the orthogonal complement
$\Pic(V(B))_0$ of $K_{V(B)}$ in $\Pic(V(B))$ coincides with
$$\{a_0 l_0+\sum_{b\in B}a_b l_b\mid a_0,a_b \in\Z, -3 a_0+\sum_{b\in B}a_b
=0\};$$ it is a Galois-invariant {\sl pure} free commutative
subgroup of rank $7$.

Notice that one may view $C_B$ as a $K$-curve on $V(B)$ \cite[p.
160]{Dolgachev}. Then the inclusion map $C_B\subset V(B)$ induced
the homomorphism of Galois modules
$$r:\Pic(V(B)) \to \Pic(C_B)$$
where $\Pic(C_B)$ is the Picard group of $C_B$ over $K_a$. Recall
that $J(B)(K_a)$ is a Galois submodule of $\Pic(C_B)$ that
consists of divisor classes of degree zero. In particular,
$J(B)_2$ coincides with the kernel $\Pic(C_B)_2$ of multiplication
by $2$ in $\Pic(C_B)$. It is known (Lemma 1 on p. 161 of
\cite{Dolgachev}) that
$$r(\Pic(V(B))_0)\subset \Pic(C_B)_2=J(B)_2.$$
This gives rise to the homomorphism
$$\bar{r}:\Pic(C_B)_0/2\Pic(C_B)_0 \to J(B)_2, \quad
D+2\Pic(C_B)_0\mapsto r(D)$$ of Galois modules. By Lemma 2 on pp.
161-162 of \cite{Dolgachev}, the kernel of $\bar{r}$ is as
follows. The intersection form on $\Pic(V(B))$ defines by
reduction modulo $2$ a symmetric bilinear form
\begin{gather*}\psi:\Pic(V(B))/2\Pic(V(B)) \times \Pic(V(B))/2\Pic(V(B))
\to\Z/2\Z=\F_2,\\ \quad D+2\Pic(V(B)) ,D'+2\Pic(V(B))\mapsto
(D,D')+2\Z\end{gather*} and we write
$$\psi_0:\Pic(V(B))_0/2\Pic(V(B))_0 \times \Pic(V(B))_0/2\Pic(V(B))_0
\to\F_2$$ for the restriction of $\psi$ to $\Pic(V(B))_0$.
 Then the kernel (radical)
of $\psi_0$ coincides with $\ker(\bar{r})$. (The same Lemma also
asserts that $\bar{r}$ is surjective.)

Let us describe explicitly  the kernel of $\psi_0$. Since
$\Pic(V(B))_0$ is a pure subgroup of $\Pic(V(B))$, we may view
$\Pic(V(B))_0/2\Pic(V(B))_0$ as a $7$-dimensional $\F_2$-vector
subspace (even Galois submodule) in $\Pic(V(B))/2\Pic(V(B))$. Let
$\bar{l}_0$ (resp. $\bar{l}_b$) be the image of $l_0$ (resp.
$l_b$) in $\Pic(V(B))/2\Pic(V(B))$. Then
$\{\bar{l}_0,\{\bar{l}_b\}_{ b\in B}\}$ constitute an {\sl
orthonormal} (with respect to $\psi$) basis of the $\F_2$-vector
space $\Pic(V(B))/2\Pic(V(B))$. Clearly, $\psi$ is non-degenerate.
It is also clear that
$$\Pic(V(B))_0/2\Pic(V(B))_0=\{a_0 \bar{l}_0+\sum_{b\in B}a_b \bar{l}_b\mid a_0,a_b \in\F_2,  a_0+\sum_{b\in B}a_b
=0\}$$ is the orthogonal complement of {\sl isotropic}
$$\bar{v}_0=\bar{l}_0+\sum_{b\in B}\bar{l}_b$$
in $\Pic(V(B))/2\Pic(V(B))$ with respect to $\psi$. Notice that
$\bar{v}_0$ is Galois-invariant. The non-degeneracy of $\psi$
implies that the kernel of $\psi_0$ is the Galois-invariant
one-dimensional $\F_2$-subspace generated by $\bar{v}_0$.

This gives us the injective homomorphism
$$(\Pic(V(B))_0/2\Pic(V(B))_0)/\F_2 \bar{v}_0 \hookrightarrow
J(B)_2$$ of Galois modules; dimension arguments imply that it is
an isomorphism. So, in order to finish the proof, it suffices to
construct a surjective homomorphism $\Pic(V(B))_0/2\Pic(V(B))_0
\twoheadrightarrow Q_B$ of Galois modules, whose kernel coincides
with $\F_2 \bar{v}_0$. In order to do that, let us consider the
homomorphism $$\kappa:\Pic(V(B))_0/2\Pic(V(B))_0 \to Q_B$$ that
sends $z=a_0 \bar{l}_0+\sum_{b\in B}a_b \bar{l}_b$ to the function
$\kappa(z):b \mapsto a_b+a_0$. Since
$$a_0+\sum_{b\in
B}a_b=0  \ \text{and}  \ \#(B)a_0 =7a_0=a_0\in\F_2,$$  indeed we
have $\kappa(z)\in Q_B$. It is also clear that $\kappa(z)$ is
identically zero if and only if $a_0=a_b \ \forall \ b$, i.e.
$z=0$ or $\bar{v}_0$. Clearly, $\kappa$ is a surjective
homomorphism of Galois modules and $\ker(\kappa)=\F_2 \bar{v}_0$.
\end{proof}

\begin{proof}[Proof of Lemma \ref{general}]
We will use a notation $(x:y:z)$ for homogeneous coordinates on
$\P^2$.  Suppose that here are three points of $B_f$ that lie on a
line $ax +by+cz=0$. This means that there are distinct roots
$\alpha_1,\alpha_2,\alpha_3$ of $f$ and elements $a,b,c \in K_a$
such that all $a \alpha_i^3+b \alpha_i+c=0$ and, at least, one of
$a,b,c$ does not vanish. It follows that the polynomial
$at^3+bt+c\in K_a[t]$ is not identically zero and has three
distinct roots $\alpha_1,\alpha_2,\alpha_3$. This implies that
 $a\ne 0$ and
$$at^3+bt+c=a(t-\alpha_1)(t-\alpha_2)(t-\alpha_3).$$ It follows that  $\alpha_1+\alpha_2+\alpha_3=0$. Let
us denote the remaining roots of $f$ by
$\alpha_4,\alpha_5,\alpha_6,\alpha_7$. Clearly, $\Gal(K)$ acts
$3$-transitively on $\RR_f$. This implies that there exists
$\sigma\in\Gal(K)$ such that
$$\sigma(\alpha_1)=\alpha_4,\sigma(\alpha_2)=\alpha_2,\sigma(\alpha_3)=\alpha_3$$
and therefore
 $\alpha_2+\alpha_3+\alpha_4=\sigma(\alpha_2+\alpha_3+\alpha_1)=0$ and
therefore $\alpha_1=\alpha_4$ which is not the case. The obtained
contradiction proves that no three points of $B_f$ lie on a  one
line.

Suppose that six points of $B_f$ lie on a one conic. Let
$$a_0 z^2+a_1 yz+a_2 y^2+a_3 xz+a_4 xy + a_6x^2 =0$$ be an equation of the conic. Then not all $a_i$ do vanish and there are six
distinct roots $\alpha_1, \cdots ,\alpha_6$ of $f$ such that all
$a_6 \alpha_k^6+\sum_{i=0}^4 a_i\alpha_k^i=0$. This implies that
the polynomial $a_6t^6+\sum_{i=0}^4 a_i t^i$ is not identically
zero and has $6$ distinct roots $\alpha_1, \cdots \alpha_6$. It
follows that $a_6 \ne 0$ and
$$a_6t^6+\sum_{i=0}^4 a_i t^i=a_6 \prod_{i=1}^6 (t-\alpha_i).$$ This implies that
$\sum_{i=1}^6 \alpha_i=0$. Since the sum of all roots of $f$ lies
in $K$,  the remaining seventh root of $f$ lies in $K$. This
contradicts to the irreducibility of $f$. The obtained
contradiction proves that no six points of $B_f$ lie on a one
conic.
\end{proof}

\begin{lem}
\label{endo1} Let $B\subset \P^2(K_a)$ be a $7$-element set of
points in general position. Assume that $\Gal(K)$ permutes
elements of $B$ and
 the image of $\Gal(K)$ in $\Perm(B)\cong \mathbf{S}_7$ coincides either
with the full symmetric group $\mathbf{S}_7$ or with with the alternating
group $\A_7$; in particular, $B$ consitutes a Galois orbit.  Then
either $\End(J(B))=\Z$ or $\fchar(K)>0$ and $J(B)$ is a
supersingular abelian variety.
\end{lem}

\begin{proof}[Proof of Lemma \ref{endo1}]
Recall that $G_B$ is the image of $\Gal(K)$ in $\Perm(B)$. By
assumption, $G_B=\mathbf{S}_7$ or $\A_7$. It is known \cite[Ex.
7.2]{ZarhinTexel} that the $G_B$-module $Q_B$ is {\sl very simple}
in the sense of \cite{ZarhinTexel,ZarhinVery,ZarhinMMJ}. In
particular,
$$\End_{G_B}(Q_B)=\F_2.$$ The surjectivity of
$\Gal(K)\twoheadrightarrow G_B$ implies that the $\Gal((K)$-module
$Q_B$ is also very simple. Applying Lemma \ref{points}, we
conclude that the $\Gal((K)$-module $J(B)_2$ is also very simple.
Now the assertion follows from lemma 2.3  of \cite{ZarhinTexel}.
\end{proof}

\begin{proof}[Proof of Theorem \ref{endo}]
In light of Lemma \ref{endo1}, we may and will assume that
$\fchar(K)>0$ and $J(B)$ is a supersingular abelian variety. We
need to arrive to a contradiction. Replacing if necessary $K$ by
its suitable quadraric extension we may and will assume that
$G_B=\A_7$. Adjoining to $K$ all $2$-power roots of unity, we may
and will assume that $K$ contains all $2$-power roots of unity and
still  $G_B=\A_7$. It follows from Lemma \ref{points} that $\A_7$
is isomorphic to the image of $\Gal(K) \to \Aut_{\F_2}(J(B)_2)$
and the $\A_7$-module $J(B)_2$ is very simple; in particular,
$\End_{\A_7}(J(B)_2)=\F_2$. Applying Theorem 3.3 of
\cite{ZarhinBSMF}, we conclude that there exists a central
extension $G_1 \twoheadrightarrow \A_7$ such that $G_1$ is
perfect, $\ker(G_1\twoheadrightarrow \A_7)$ is a central cyclic
subgroup of order $1$ or $2$ and there exists a symplectic
absolutely irreducible $6$-dimensional representation of $G_1$ in
characteristic zero. This implies (in notations of \cite{Atlas})
that either $G_1\cong \A_7$ or $G_1\cong 2.\A_7$. However,  the
table of characters on p. 10 of \cite{Atlas} tells us that neither
$\A_7$ nor $2.\A_7$ admits a {\sl symplectic} absolutely
irreducible $6$-dimensional representation in characteristic zero.
The obtained contradiction proves the Theorem.
\end{proof}

\section{Explicit formulas}
In this section we describe explicitly $H_B$ when $B=B_f$. We have
$$f(t)=\sum_{i=0}^7 c_i t^i \in K[t], \ c_7 \ne 0.$$
 We are going to
describe explicitly  cubic forms that vanish on $B_f$. Clearly,
$u:=xz^2-y^3$ and $v:=c_7 x^2 y+ c_6 x^2 z+ c_5 x y^2+c_4 xyz+c_3
x z^2+ c_2 y^2 z +c_1 y z^2 + c_0 z^3$ vanish on $B_f$. In order
to find a third vanishing cubic form, let us define a polynomial
$h(t)\in K[t]$ as a (non-zero) remainder with respect to division
 by $f(t)$:
$$t^9 -h(t) \in f(t) K[t],\ \deg(h)<\deg(f)=7.$$
We have
$$h(t)=\sum_{i=0}^6 d_i t^i\in K[t].$$
For all roots $\alpha$ of $f$ we have
$$0=\alpha^9-h(\alpha)=\alpha^9-\sum_{i=0}^6 d_i \alpha^i.$$
This implies that the cubic form $w=x^3-d_6 x^2 z-d_5 x
y^2-d_4xyz-d_3 x z^2- d_2 y^2 z -d_1 y z^2 - d_0 z^3$ vanishes on
$B_f$. Since $u,v,w$ have $x$-degree 1,2,3  respectively, they are
linearly independent over $K_a$ and therefore constitute a basis
of $3$-dimensional $H_{B_f}$.

Now assume (till the end of this Section) that $\fchar(K)\ne
3$.\footnote{This condition was inadvertently omitted in the
Russian version \cite{ZarhinVR}.} Since $C_{B_f}$ is the
ramification curve for $\pi$, it follows that
$$g_B(C_{B_f})=\left\{(x:y:z), \begin{vmatrix}
u_x & u_y & u_z\\
v_x & v_y & v_z\\
w_x & w_y & w_z
\end{vmatrix}=0\right\}\subset \P^2$$
is a singular sextic which is $K$-birationally isomorphic to
$C_{B_f}$. (See also \cite[proposition 2 on p. 167]{Dolgachev}.)

\section{Another proof}

The aim of this Section is to give a more elementary proof of
Theorem \ref{endo} that formally does not refer to Lemma 2 of
\cite[Lemma 2 on pp. 161--162]{Dolgachev} (and therefore does not
make use of the Smith theory. However, our arguments are based on
ideas of \cite[Ch. IX]{Dolgachev}.) In order to do that, we just
need to prove Lemma \ref{points} under an additional assumption
that the image of $\Gal(K)$ in $\Perm(B)$ is ``very big".

\begin{lem}
\label{pointsB} Let $J(B)_2$ be the kernel of multiplication by
$2$ in $J(B)(K_a)$. Suppose that  $G_B$ coincides either with
$\Perm(B)$ or with $\A_7$. Then the Galois modules $J(B)_2$ and
$Q_B$ are  isomorphic.
\end{lem}

\begin{proof}
Let $g_0: V(B) \to V(B)$ be the  Geiser involution \cite[p. 66--
67]{Demazure}, i.e., the biregular covering transformation of
$\pi$. Clearly, $g_0$ is defined over $K$. This implies that if
$E$ is an irreducible $K_a$-curve on $V(B)$ then $E$ and $g_0(E)$
have the same stabilizers in $\Gal(K)$.
 Clearly, different points $b_1$ and $b_2$ of $B$ have different
stabilizers in $G_B$ and therefore in $\Gal(K)$. This implies that
$g_0(E_{b_1}) \ne E_{b_2}$, since the stabilizers of
$g_0(E_{b_1})$ and $E_{b_2}$ coincide with the stabilizers of
$b_1$ and $b_2$ respectively. This implies that the lines
$$\pi(E_{b_1}), \pi(E_{b_2}) \subset \P^2,$$
which are bitangents to $C_B$ \cite[p. 68]{Demazure}, do not
coincide.

For each $b\in B$ we write $D_b$ for the effective degree $2$
divisor on the plane quartic $C_B$ such that $2 D_b$ coincides
with the intersection of $C_B$ and $\pi(E_{b})$; it is well known
that (the linear equivalence class of) $D_b$ is a theta
characteristic on $C_B$. It is also clear that
$$\sigma(D_b)=D_{\sigma(b)} \quad \forall\ \sigma\in\Gal(K),\ b\in B.$$
Clearly, if $b_1\ne b_2$ then $D_{b_1}\ne D_{b_2}$ and the
divisors $2D_{b_1}$ and $2D_{b_2}$ are linearly equivalent. On the
other hand, $D_{b_1}$ and $D_{b_2}$ are {\sl not} linearly
equivalent. Indeed, if $D_{b_1}-D_{b_2}$ is the divisor of a
rational function $s$ then $s$ is a non-constant rational function
on $C_B$ with, at most, two poles. This implies that either $C_B$
is either a rational (if $s$ has exactly one pole) or
hyperelliptic (if $s$ has exactly two poles). Since a smooth plane
quartic is neither rational nor hyperelliptic, $D_{b_1}-D_{b_2}$
is  not a principal divisor.

Let $(\Z^B)^0$ be the free commutative group of all functions
$\phi:B \to\Z$ with $\sum_{b\in B}\phi(b)=0$. Clearly, $(\Z^B)^0$
is provided with the natural structure of $\Gal(K)$-module and
there is a natural isomorphism of $\Gal(K)$-modules
$$(\Z^B)^0/2 (\Z^B)^0 \cong Q_B.$$
 Let us consider the homomorphism of commutative groups
 ${\mathfrak r}: (\Z^B)^0 \to \Pic(C_B)$
that sends a function $\phi$ to the linear equivalence class of
$\sum_{b\in B} \phi(b)D_b$. Clearly, $${\mathfrak
r}((\Z^B)^0)\subset J(B)_2\subset \Pic(B)$$ and therefore
${\mathfrak r}$ kills $2\cdot (\Z^B)^0$. On the other hand, the
image of ${\mathfrak r}$ contains the (non-zero) linear
equivalence class of $D_{b_1}-D_{b_2}$. This implies that
${\mathfrak r}$ is not identically zero and we get a non-zero
homomorphism of $\Gal(K)$-modules
$$\bar{\mathfrak r}: Q_B \cong (\Z^B)^0/2 (\Z^B)^0 \to J(B)_2.$$
It is well-known that our assumptions on $G_B$ imply that the
$G_B$-module $Q_B$ is (absolutely) simple and therefore $Q_B$,
viewed as Galois module, is also simple. This implies that
$\bar{\mathfrak r}$ is injective. Since the $\F_2$-dimensions of
both $Q_B$ and $J(B)_2$ equal to $6$ and therefore coincide, we
conclude that $\bar{\mathfrak r}$ is an isomorphism.
\end{proof}

\section{Added in translation}
The following assertion is a natural generalization of Lemma
\ref{general}.

\begin{prop}\label{General}
Suppose that $E\subset \P^2$ is an absolutely irreducible cubic
curve that is defined over $K$. Suppose that $B\subset E(K_a)$ is
a a $7$-element set that is a $\Gal(K)$-orbit. Let us assume that
 the image $G_B$ of $\Gal(K)$ in the group $\Perm(B)$
of all permutations of $B$ coincides either with
$\Perm(B)\cong\mathbf{S}_7$ or with the alternating group $\A_7$. Then
$B$ is in general position.
\end{prop}

\begin{proof}
Clearly, $\Gal(K)$ acts $3$-transitively on $B$.

 Step 1. Suppose
that $D$ is a line in $\P^2$  that contains three points of $B$
say,
$$\{P_1,P_2,P_3\}\subset \{P_1,P_2,P_3, P_4,P_5,P_6,P_7\}=B.$$
Clearly, $D\bigcap E=\{P_1,P_2,P_3\}$. There exists $\sigma
\in\Gal(K)$ such that $\sigma(\{P_1,P_2,P_3\})=\{P_1,P_2,P_4\}$.
It follows that the line $\sigma(D)$ contains $\{P_1,P_2,P_4\}$
and therefore $\sigma(D)\bigcap E=\{P_1,P_2,P_4\}$. In particular,
$\sigma(D)\ne D$. However, the distinct lines $D$ and $\sigma(D)$
meet each other at {\sl two} distinct points $P_1$ and $P_2$.
Contradiction.

Step 2. Suppose that $Y$ is a conic in $\P^2$ such that $Y$
contains six points of $B$ say, $\{P_1,P_2,P_3,
P_4,P_5,P_6\}=B\setminus\{P_7\}$. Clearly, $Y\bigcap
E=B\setminus\{P_7\}$. If $Y$ is reducible, i.e., is a union  of
two lines $D_1$ and $D_2$ then either $D_1$ or $D_2$ contains (at
least) three points of $B$, which is not the case, thanks to Step
1. Therefore $Y$ is {\sl irreducible}.

There exists $\sigma \in\Gal(K)$ such that $\sigma(P_1)=P_7$. Then
$\sigma(P_7)=P_i $ for some positive integer $i\le 6$. This
implies that $\sigma(B\setminus\{P_7\})=B\setminus\{P_i\}$ and
 the irreducible conic $\sigma(Y)$ contains
$B\setminus\{P_i\}$. Clearly, $\sigma(Y)\bigcap
E=B\setminus\{P_i\}$ contains $P_7$. In particular, $\sigma(Y)\ne
Y$. However, both conics contain the $5$-element set
$B\setminus\{P_i,P_7\}$. Contradiction.
\end{proof}


\begin{thebibliography}{99}

\bibitem{Atlas} J. H. Conway, R. T. Curtis, S. P. Norton, R. A. Parker, R. A. Wilson, Atlas of finite groups. Clarendon Press, Oxford, 1985.

\bibitem{Demazure} M. Demazure, {\em Surfaces de Del Pezzo} II,III,IV,V.
Springer Lecture Notes in Math. {\bf 777} (1980), 23--69.

\bibitem{Dolgachev} I. Dolgachev, D. Ortland, Point sets in
projective spaces and theta functions, Ast\'erisque {\bf 165}
(1986).

\bibitem{Dolgachev2}\footnote{This reference was updated.}
 I. Dolgachev,   Classical Algebraic geometry. A modern view.
Cambridge University Press, 2012. (Section 8.7)

\bibitem{Har} R. Hartshorne, Algebraic Geometry, GTM {\bf 52},
Springer-Verlag, 1977.

\bibitem{Iskovskikh} V. A. Iskovskikh, {\em Minimal models of
rational surfaces over arbitrary fields}. Izv. Akad. Nauk Ser.
Mat. {\bf 43} (1979), 19--43; Math. USSR-Izv.  {\bf 14} (1980),
17--39.

\bibitem{IskovskikhSh} V. A. Iskovskikh, I.R. Shafarevich, {\em
Algebraic surfaces}.  Algebraic geometry, II,  127--262,
Encyclopaedia Math. Sci., {\bf 35}, Springer, Berlin, 1996.

\bibitem{Manin} Yu. I. Manin, Cubic forms, second edition, North
Holland, 1986.

\bibitem{T} Yu. I. Manin,  M. A. Tsfasman,
{\em Rational varieties: algebra, geometry, arithmetic}.  Uspekhi
Mat. Nauk {\bf 41} (1986), no. 2(248), 43--94; Russian Math.
Surveys {\bf 41} (1986), no. 2, 51--116.

\bibitem{ZarhinMRL} Yu. G. Zarhin, {\em Hyperelliptic jacobians without complex multiplication}.
 Math. Res. Letters {\bf 7} (2000), 123--132.

\bibitem{ZarhinTexel} Yu. G. Zarhin, {\em Hyperelliptic jacobians and modular
representations}. In: Moduli of abelian varieties (eds. C. Faber,
G. van der Geer and F. Oort). Progress in Math., vol. {\bf 195}
(Birkh\"auser, 2001), pp. 473--490.

\bibitem{ZarhinMRL2}  Yu. G. Zarhin, {\em Hyperelliptic jacobians without
complex multiplication in positive characteristic}. Math. Res.
Letters {\bf 8} (2001), 429--435.



\bibitem{ZarhinMMJ} Yu. G. Zarhin, {\em Very simple $2$-adic representations
 and hyperelliptic jacobians}. Moscow Math. J.
  {\bf 2} (2002), issue 2, 403-431.


\bibitem{ZarhinVery} Yu. G. Zarhin, {\em Very simple representations: variations on a theme of
Clifford}. In:  Progress in Galois Theory (H. V\"olklein, T.
Shaska eds.), Springer Verlag, 2005, pp. 151--168.



\bibitem{ZarhinBSMF} Yu. G. Zarhin, {\em Non-supersingular hyperelliptic
jacobians}. Bull. Soc. Math. France {\bf 132} (2004), 617--634.

\bibitem{ZarhinVR} Yu. G. Zarhin, {\em Del Pezzo surfaces of degree $2$ and  jacobians without complex
multiplication} (Russian). Trudy St. Petersburg  Mat. Obsch. {\bf
11} (2005), 81--91.
\end{thebibliography}
\end{document}